\documentclass[12pt]{amsart}

\usepackage{graphicx}
\usepackage{amssymb}
\usepackage{setspace}
\usepackage{amsmath}     
\usepackage{stackrel}    
\usepackage{stackengine} 

\usepackage[all]{xy}

\usepackage{mathtools}
\usepackage{amsfonts}

\usepackage{tikz-cd}

\newcommand{\C}[1]{{\mathcal#1}} 

\theoremstyle{plain}

\theoremstyle{definition}

\theoremstyle{definition}

\theoremstyle{remark}

\usepackage{fancyhdr}
\begin{document}
\footskip=30pts

\title{On p-adic Frobenius lifts and p-adic periods,\\ from a Deformation Theory viewpoint} %
\author{Lucian M. Ionescu}
\address{Department of Mathematics, Illinois State University, IL 61790-4520}
\email{lmiones@ilstu.edu}
\date{\today} 

\begin{abstract}
Presenting p-adic numbers as {\em deformations} of finite fields
allows a better understanding of Frobenius lifts
and their connection with p-derivations in the sense of Buium \cite{Buium-Main}.
In this way ``numbers {\em are} functions'', as recognized before \cite{Manin:Numbers},
allowing to view initial structure deformation problems
as arithmetic differential equations as in \cite{Buium-Manin}, and 
providing a cohomological interpretation to Buium calculus via Hochschild cohomology
which controls deformations of algebraic structures.

Applications to p-adic periods are considered, 
including to the classical Euler gamma and beta functions and their p-adic analogues,
from a cohomological point of view.

Connections between various methods for computing scattering amplitudes are 
related to the moduli space problem and period domains.
\end{abstract}

\maketitle
\setcounter{tocdepth}{3} 
\tableofcontents

\section{Introduction}  
Algebraic integrals called periods \cite{K,Muller-Stach:Period,Carlson:PD}\footnote{See \cite{KZ} for additional aspects regarding periods.}, 
are a new class of numbers extending algebraic numbers and bridging quantum physics and mathematics together, in a way yet to be better understood \cite{Brown:Periods-Feynman,Schnetz:QuantumPeriods,Ionescu-Sumitro}.

A p-adic version for periods has been defined, in quite abstract terms, in connection with Tate-Hodge Conjecture and later on, Fontaine Conjectures \cite{Wiki:p-adic Hodge Theory}, but a more ``down-to-earth'' lucrative approach is yet to be found, in order to support the ``movement'' of physics towards an ``Ultimate Physics Theory'' in terms of Number Theory, as for example the developments of a {\em p-adic String Theory}
\cite{Volovich:NT,Dragovich:p-adic Math}, 
and one which ``sounds'' more familiar to the physicist, who expects the keyword ``quantization'' somewhere down the path (and ``deformation'' can be snicked quite nicely before that!).

From a computer science perspective, instead, a ``Quantum Theory'' could be formulated already ``quantum'', i.e. {\em discrete} and of finite type, using the language of {\em Quantum Computing} at foundations, as a general {\em flow on networks} calculus (e.g. Turaev's on Ribbon Categories), at a conceptual level
\footnote{Indeed, what are {\em quantization functors} on cobordism categories, or spin networks etc. if not that!?}.
What still is needed is to understand the role of primes (and Riemann spectrum),
in the context of what the present author calls {\em p-adic (finite) strings}.
But then, a {\em deformation quantization} of the theory of finite fields should be beneficial for a brainstorming translation of the ``good concepts'' from ``continuous mathematical-physics over complex numbers, to the discrete graded case of p-adics and adeles'.
Indeed, if the complex numbers are thought of as the algebraic closure of the topological closure of the rationals $C=\bar{Q}_\infty$, then algebraic extensions of p-adic numbers should suffices $Q_q, q=p^n, n\in N$, even though they are not algebraically closed (the closure $C_p$ is much bigger: see \cite{Koblitz}).

In this note we will explain how p-adic numbers can be understood as deformations of finite fields as Galois-Klein geometries, in terms close to the spirit of deformation quantization, a perspective hopefully beneficial for incorporating Algebraic Geometry into a Number Theory approach to {\em Finite String Theory}.

The article is organized as follows. 
The p-adic numbers are first reformulated as formal series, in the sense of Deformation Theory, together with basic facts, including how a canonical lift of Frobenius, as a companion deformation of their symmetries \cite{Wiki:Galois Deformations}. 

Then the p-derivations forming the basis of the lambda calculus of Buium are reviewed, and reinterpreted, from the deformation point of view. The appearance of Hochschild cohomology is no accident,
and allows to address the periods and the period comparison isomorphism in a new light, briefly.

As applications, the article comments on the p-adic analogues of Euler's gamma and beta integrals, which are natural extensions of the Gauss and Jacobi sums in finite characteristic.
The connections with scatering amplitudes, as Feynman Integrals and Veneziano amplitudes, well know to be ``coincidentaly'' related to number theoretic values like Multiple Zeta Values, no longer seem so unexpected, in view of the correspondence between the theories of period isomorphisms.

In conclusion, the author ponders on the similarity between the familiar mathematics of real numbers and that of p-adics ``analysis'', or rather deformation theory, which substantiates the claim that p-adic numbers {\em are} in fact functions, will all their benefits.
Finally, further directions of development are discussed, and questions are raised, especially regarding the {\em moduli problem} from a Deformation Theory standpoint.

\section{p-Adic Numbers as Functions}
The traditional Cauchy completion approach to p-adic numbers $Q_p$ hides their alternative algebraic presentation as a $Z$-module of formal series with coefficients in $(Z/pZ,+)$, as noted long time ago by Hensel, who compared them with Laurent series:
$$Z_p=(F_p[[h]],*)\ni a=\sum_{i=1}^\infty a_nh^n, \quad a_n \in F_p, 
\quad a_n*b_n=(a_n+_p b_n) + c_p(a_n,b_n)h.$$
We will focus on the ``integral elements'', since the filed extensions occur as corresponding fraction fields\footnote{There is a ``complication'' due to a different convenient definition of integral closure, which is of no concern here.}.

Note also that the field structure $(F_p,+,\cdot)$ on $(Z/pZ,+_p)$ comes for free as ``repeated addition'', and that it corresponds to the group of automorphisms as an abelian group (``discrete space of arrows-vectors''). 

\subsection{p-adic Deformations of Finite Fields}
The deformation parameter is the generic ``Planck's constant'' $h$ to avoid the miss-conception that they are ``just'' the integers completed in the ``wrong direction'', and to emphasize the conceptual role of the ``uniformizer'' as a grading parameter.
On the contrary, the real numbers result in the {\em topological completion} in the wrong direction, of ``infinity-small'', not needed in a quantum world of atoms and Zeno paradoxes \cite{RealFish}.
Note also that the completion in the direction of the carry-over 2-cocyle 
$$c_p(x,y):F_p\times F_p\to F_p,\
c_p(x,y)=[j(x)+j(y)\ mod_p]-(x+_p y),\ mod_p:Z\leftrightarrow F_p: j,$$
which has an awkward ``closed form'' in terms of a section $j$, as one already knows from the need of introducing Witt vectors \cite{Finkel-Witt}.

To be more specific, we ``zoom in'' on the infinitesimal deformation for now, i.e. let us consider the infinitesimal deformation of $Z/pZ$, or equivalently its central extension with 2-cocycle $c_p$ \cite{Gerstenhaber}:
$$0 \to Z/p \to Z/p^2\to Z/p\to 0.$$
Then the addition modulo $p$ is deformed into a ``$*$-product'' as above, to use the physics jargon when comes to quantization via deformation. 
Indeed, at ``tangent space'' level (degree 0), the central extension is defined as follows:
$$(a_0,0)*(b_0,0)=(a_0+b_0, c_p(a_0,b_0)),$$
where we have dropped the subscript from $+_p$.

\subsection{The Canonical Frobenius Lift}
Having deformed the abelian group $(Z/pZ,+)$, 
apply a general principle, that an object is determined by its symmetries.
Indeed, the Frobenius automorophisms generates the Galois extension $[F_{p^n}:F_p]$,
so, the deformed structure $Z_q=F_{p^n}[[h]]$ (with $q=p^n$), will have automorphisms generated by the corresponding deformation \cite{Wiki:Galois Deformations}.

In order to have a nice commutation relation, we need to use the Teichmuler character
to represent the p-adic digits, i.e. represent the formal deformations as 
Witt vectors \cite{Rabinoff}, p.14:
$$Z_q\cong W(F_q) \ni x=\sum \tau(x_n)p^n.$$
Then the {\em canonical lift of Frobenius} $\phi$ ``commutes'' with the Teichmuler character:
$$\phi_p(\tau(v))=\tau(v)^n, \ v\in F_q, \quad \phi(x)=\sum \tau(x_n)^n p^n,\ x=\sum \tau(x_n)p^n,$$
i.e. the Teichmuler character pull-back of the lift of Frobenius is the Frobenius pull-back of the lift of the Teichmuler character.


\section{Relations with Buium Calculus}
As mentioned above, p-adic numbers are analogues of Laurent series, as noticed initially by Hansel, and more recently by Manin \cite{Manin:Numbers},
idea developed recently into a full-blown calculus by Buium \cite{Buium-Main}.

\subsection{p-Derivations and Lambda Calculus}
Following \cite{Buium-Manin} (specialixed to p-adics framework, 
with $f(x)=x$ be identity\footnote{Extensive details can be found in \cite{Buium-Main}.}, 
a {\em p-derivation} $\delta(x)$ on p-adic numbers, thought of as functions conform \cite{Manin:Numbers},
is a mapping satisfying the following conditions \cite{Buium-Manin}, p.2:
$$\delta(x+y)=\delta_p(x)+\delta_p(y)+C_p(x,y), \quad C_p(x,y)=[x^p+y^p-(x+y)^p]/p,$$
$$\delta(xy)=x^p\delta(y)+\delta_p(x)y^p+ p\cdot \delta_p(x)\delta_p(y).$$

It is immediate that for any such p-derivation allows to define a 
ring endomorphism $\phi_p(x)=x^p+p \delta_p(x)$,
which is a {\em lift of Frobenius} 
$\phi_p(x)=x^p\ mod \ p$ (loc. cit.).

\subsection{Relation with Deformation Theory}
With our notation, $p\mapsto h$ and $+\mapsto *$, 
the cohomological interpretation of the 1st equation above is
$d^*_{Hoch}\delta_p=c_p$, where the associative operation in view is p-adic addition
$*$, i.e. the deformation of component-wise addition of formal power series in $h$,
while the 2nd equation is, similarly,
$d^\cdot_{Hoch}\delta_p=0$, when viewing $Z_q$ as an $Z_q$-module via the ``natural'' Frobenius action $L_x(y)=x^ny$. 
Now rewrite a p-derivation using our notation, as a Kahler differential 1-form $\omega(x)=\delta(x)dh$, in the formal series deformation parameter $h$:
$$\omega=[\phi_p(x)-Frob(x)]\ dh/h.$$
Then it is a {\em Cauchy kernel}, i.e. the product of an ``entire function'' times the generator $dh/h$ of the 1-st (Monski-Washnitzer) algebraic de Rham homology group $H^1_{MW}(X)$ (see \cite{Hartog}, Ex. 3.1.9, p.28).

\vspace{.2in}
The proper reinterpretation of the above formulas in terms of the Teichmuler character and Frobenius lift will be discussed separately \cite{LI:Periods-FI-JS}.

\section{Applications to Periods}
Nevertheless the deformation theory approach sketched above allows for a better understanding of the p-adic periods, in view of a comparison of the MW-cohomology with Hochschild cohomology \cite{Anel:Deformation-Cohomology}, 
and corresponding period isomorphism \cite{Hartog,Kedlaya}. 

\subsection{p-adic Euler functions analogues}
The p-adic analogue of Euler's Gamma, 
the so called Morita Gamma function, and of the Beta function are \cite{Boyarsky}:
$$\Gamma_p= (-1)^x \prod\limits_{1<i<x, p\not| i} i, \qquad B_p(x)=\Gamma_p(a)\Gamma_p(b)/\Gamma(a+b)$$
having the following Hochschild cohomological interpretation $Beta=d_H\Gamma$.

If an analogue of the p-adic Gauss sum is defined (Eq. (4.1) loc. cit. p.362, with $f=1$ for simplicity)
$g_p=\tau \circ \theta$, as a composition of Teichmuler character $\tau$
and the exponential $\theta(x)=exp[\pi (x-x^q)]$, where $\pi^{p-1}=-1$,
then a similar cohomological interpretation of a 2-cocycle holds: $B_p=d\Gamma_p$ (p.366).

Before commenting on the role of periods of some special values of these functions,
we further note that a p-adic analog of the Jacobi sum, naturally deformed via the Teichmuler character (compare with Eq. (7.2) \cite{Boyarsky}, p365, except for the chosen negative sign):
$$J(c,c')=(c*c')(1), \quad c=\tau^\alpha, c'=\tau^\beta,$$
is again, in a consistent way, an exact 2-coboundary 
$J_p=d_{Hoch}g_p$\footnote{Signifies the absence of obstructions 
for a full formal deformation.}.

The author's interpretation in the context of Monsky-Weishnizter cohomology,
an algebraic de Rham adaptation to the p-adic case,
may be understood in the context of period isomorphisms, next.

\subsection{Relation to Period Isomorphisms}
The de Rham period isomorphism for algebraic varieties over complex numbers have various analogues to p-adics, as developed by Tate, Faltings and more recently Fontaine \cite{Wiki:p-adic Hodge Theory}, who introduced the so call big $B$-rings of periods \cite{Fontaine}.
The deformation theory viewpoint adopted above suggests a possible comparison isomorphim can be formulated in terms of Hochschild cohomology, instead of, and corresponding to the etale cohomology comparison isomorphism as usual.
A starting point, only, would be a cohomological interpretation of Rohrlich formula for the periods of Fermat curves $F_m:x^m+y^m=1$ in the complex case, as a product of of a Hochschild 2-coboundary and a cyclotomic number \cite{Kashio-Fermat periods}, \S 2, p.2:
$$\int_\gamma \eta_{(r/m,s/m)}\in B(r/m,s/m) \cdot Q(\zeta_m), $$
$$ \gamma \in H_1(F_m(C),Q), \quad \eta_{(r/m,s/m)}=dx^r/y^{m-s}\in H_{dR}(F_m(C),Q),$$
where $s,t\in Z/mZ$, and $s+t=1$ to suggest the underlying convolution, as in a Jocobi sum or Euler beta function, with cohomological significance.

The p-adic analog for formula for the p-adic periods allows to infer that the image of the p-adic comparison cohomology pairing has a smaller image that the ``big'' Fontaine B-rings (loc. cit. p.3):
$$H_1(F_m(C),Q)\times H^1_{dR}(F_m,Q)\to B_{dR}.$$
It is interesting to ponder on a possible interpretation of the p-adic periods derived from the action of the Frobenius lift on cohomology, conform to Lefschetz Fixed Point Theorem \cite{Hartog}, p.28:
$$\#\bar{X}(F_q)=q^{n-1}\sum_0^{n-1}(-1)^Tr\left((\phi^*)^{-1}|H^i|{MW}(\bar{A})\right),
$$
yielding a possible connection between its eigenvalues, the ``Weil zeros'' satisfying the p-adic analog of Riemann Hypothesis, and p-adic periods.

\subsection{Relation to Veneziano Amplitude}
The formal analogy with the Jacobi sum may be further strengthen by considering the
above p-adic Jacoby sum $J_p=d_{Hoch}g_p$, a 2-characters convolution value.

The Veneziano's ``educated guess'' for a Regge trajectory compatible scattering amplitude \cite{Veneziano Amplitude,Kholodenko:VA}, which in some sense ``started'' the String Theory movement, is an iterated integral on moduli spaces of punctured Riemann spheres $\C{M}_{0,4}$ \cite{Brown-Moduli}, Eq. (1.1), p.1, , with $n=1$:
$$ A(4-points)=B(a,b)=(d_{Hoch}\Gamma)(a,b)=\int\limits_{0\le t\le 1}t^a (1-t)^b dt,$$
and hence a period.
This can be related to Multiple Zeta Values, explaining in an indirect way the ``coincidence'' with the Feynman amplitudes, as linear combinations of MZVs \cite{QuantaMagazine}, 
``closing'', in a way, this circle of ideas.

It suggests that perhaps an analogous moduli space for p-adic curves exists (beyond the scope of this article), or at least a direct connection with the Multiple Zeta Values \cite{Brown-MZV}, and opening a research direction for explaining the ``unreasonable'' effectiveness of the new methods for computing scattering amplitudes, starting with the BCFW-recursion method \cite{BCFW}, towards the general, yet elusive as for now, amplituhedron concept \cite{Nima,Brown:Periods-Feynman,LI:Periods-FI-JS}.

\section{Conclusions and Further Developments}
As an overarching theme, quantization was always thought of as a sort of ``deformation'' of Newtonian physics in the ``direction'' of Planck constant, not to mention the other deformation, of Galilean Relativity in the direction of $1/c$, the inverse of the speed of light.

That {\em Deformation Theory} is a natural generalization of Lie Theory from the framework of Lie algebras / Lie groups to quite general algebraic structure \cite{LI:Lie-Def}, including the modern mathematics of Quantum Groups, 
should be enough incentive for ``prefering'' to advertise p-adic ``analysis'' as {\em p-adic deformation theory}. 

The interplay between algebraic and analytic, pertaining to field extensions $F_q$ and h-adic completions $Z_q$, with their fields of fractions $Q_q$\footnote{... perhaps better projectively completed as $Z_qP^1$ for Algebraic Geometry treatments.}, 
can be pictured as a 2-dimensional 1st octant grid of deformations of ``vector bundles'':
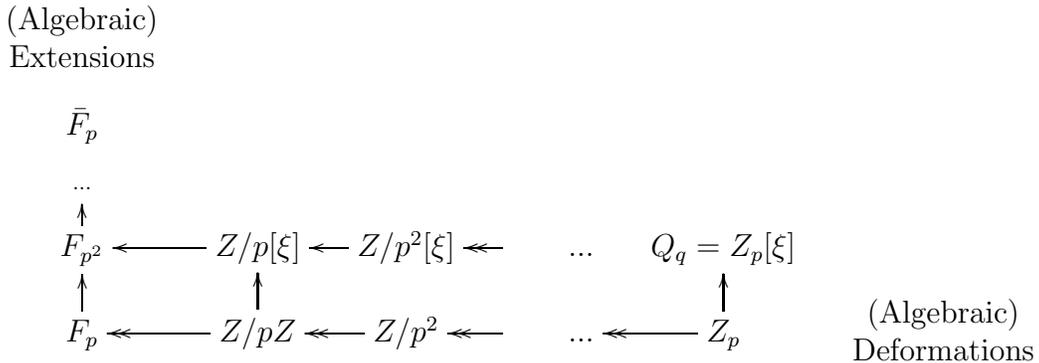
$$\xymatrix@R=.1in @C=.2in{
\stackanchor{(Algebraic)}{Extensions} & \\
\bar{F}_p & \\
\overset{...}{\ } & \\
F_{p^2} \ar[u] & \ar[l] Z/p[\xi] & \ar[l] Z/p^2[\xi] & \ar@{->>}[l] & ... & Q_q=Z_p[\xi]\\
F_p \ar[u] & \ar[u] \ar@{->>}[l] Z/pZ & \ar@{->>}[l] Z/p^2 & \ar@{->>}[l] & ... & \ar@{->>}[l] \ar[u] Z_p & \stackanchor{(Algebraic)}{Deformations}
}$$

Then the similarity between ``real/complex (topologically completing the rationals at the place/prime $p=\infty$) vs. p-adic mathematics is obvious:
$$
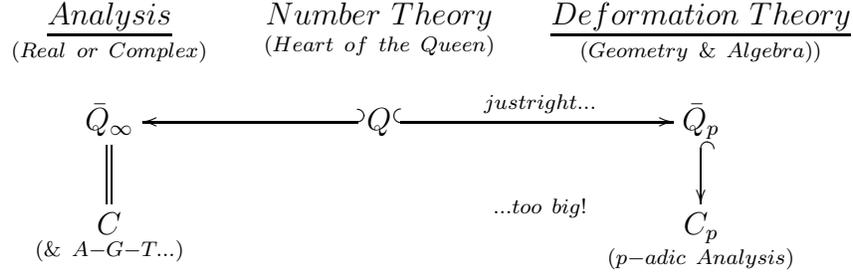
\xymatrix@R=.2in @C=.2in{
\underset{(Real \ or \ Complex)}{\underline{Analysis}} & 
\underset{(Heart\ of \ the \ Queen)}{Number\ Theory} & 
\underset{(Geometry\ \& \ Algebra))}{\underline{Deformation\ Theory}}\\
\bar{Q}_\infty \ar@{=}[d] & \ar@{_{(}->}[l] Q \ar@{^{(}->}^{just right ...}[r] & \bar{Q}_p \ar@{^{(}->}[d] \\
\underset{(\&\ A-G-T ...)}{C} & \ar@{}[l]_{} \ar@{}[r]^{... too\ big!} & \underset{(p-adic\ Analysis)}{C_p}\\
}$$
and the ``Hamletian'' question ``To use graded or non-graded structures''? has an obvious answer, since graded structures often come with benefits, e.g. the antipode of a graded bialgebra.

The advantage of the p-adics side, of being {\em graded mathematics}, is that {\em numbers can be treated as functions}: the historical connection with the integers may be broken, and $p-adic$ ``numbers'' are just $h-adic$ Laurent series, so that the corresponding fields $Q_q$ may be treated both analytically and as number fields. 
The fact that they are not algebraically complete is not an issue, since extensions should be rather viewed as objects of a category, rather than as a big huge one-object of study.

On the concrete side, the deformation theory point of view allows to proficiently make use of Hochschild cohomology, as the one that controls deformations. 
The various p-adic analogues of special functions, starting with the Gamma Function $\Gamma$, which in fact perhaps should be thought of as the Melin transform of the exponential, in the context of Fourier Duality
\footnote{The ``paradigm'' here is pairing additive and multiplicative characters, i.e. harmonic analysis of symmetries.}, and the derived 2-cocycle, the Euler Beta function, have p-adic analogues, which seem to have deeper significance in this context.

The periods from Quantum Physics, via Feynman Integrals, MZV, have p-adic analogues if we take Veneziano Amplitude as a model, and recall that it is an iterated integral on an algebraic variety with a divisor (Riemann sphere with marked points).
Than the ``translation'' of the theory of periods in p-adic realm is perhaps easier, avoiding the big $B$-period rings of Fontaine, of course, by non-experts.

Returning to the underlying theme, that ``Deformation Theory is the study of {\em infinitesimal conditions} associated to varying a solution $P$ of a problem ...'' \cite{Wiki:Deformation Theory}, explains in a way why the perturbative approach via Feynman diagrams and integrals should yield periods, and that moduli spaces are, again in some sense, {\em period domains parameterizing deformations controlled by cohomologies}.
This justifies why periods coming in families usually satisfy differential equations even in the arithmetic setup, as for example the one studied in \cite{Buium-Manin}

Also, since deformations are classified by some cohomology groups, and in particular algebraic deformations like the p-adic numbers and their endomorphisms, are characterized by Hochschild 2-cocycles with two deformations isomorphic if cohomologous \cite{Anel:Deformation-Cohomology}, p.10, Frobenius lifts as deformations of the corresponding symmetries are controlled by 1-cohomology classes determined by derivations, like the p-derivations of Buium calculus.

On a final note, understanding the adeles, i.e. how the p-adics deformations glue altogether, on the other hand, is a different story, yet Deformation Theory related, since it is perhaps about the duality of the {\em algebraic quantum group} of rationals \cite{Van Daele}.

\section{Acknowledgements}
I would like to express my gratitude for the excellent research conditions at {\em l'Institut des Hautes Etudes Scientifiques}, and for the help from Illinois State University Department of Mathematics, especially Prof. Gail Yamskulna, which made this trip possible.

The discussions with Prof. Kontsevich had a most needed guiding role, clarifying many questions of the present author, and providing extensive additional material, which only partially was used in this article; the ``straying away'' from the shown path is entirely due to author's childish attempts to explore on his own this new marvelous territory were Mathematics and Physics become one.

The encouragements from Prof. Yuri Manin, and the references provided (and some yet to be followed), confirming that, yes, ``numbers are functions'', played an important role in pursuing this line of thought.

Finally, the author also needs to thank some old friends, particularly Prof. Yan Soibelman for the lessons on Deformation Theory, and of course, Professors Louis Crane and David Yetter for TQFTs and all that jazz, knotted and braided.


\end{document}